\documentclass[a4paper,12pt,reqno]{amsart}
\usepackage{rotating,pdflscape,color,amssymb,subfigure,psfrag,amsmath,eufrak,bbm,epsfig}
\usepackage{a4wide,amscd}
\definecolor{vdarkred}{rgb}{0.6,0,0.2}
\definecolor{vdarkblue}{rgb}{0,0.2,0.6}
\usepackage[pdftex, colorlinks, linkcolor=vdarkblue,citecolor=vdarkred]{hyperref}

\usepackage[small]{caption}

\newcommand{\Ga}{\Gamma}
\newcommand{\ga}{\gamma}

\newcommand{\F}{\mc{F}}

\newcommand{\ld}{\ldots}
\newcommand{\beg}{\begin}
\newcommand{\en}{\end}

\newcommand{\trm}{\textrm}

\newcommand{\bgt}{\begin{itemize}}
\newcommand{\ent}{\end{itemize}}
\newcommand{\ite}{\item}

\newcommand{\op}{\operatorname}
\newcommand{\eqre}{\eqref}
\newcommand{\re}{\ref}
\newcommand{\la}{\label}

\newcommand{\si}{\sigma}

\newcommand{\diag}{\operatorname{diag}}

\newcommand{\ds}{\displaystyle}

\newcommand{\p}{\mathbb{P}}

\newcommand{\Tr}{\operatorname{Tr}}

\newcommand{\Ninf}{\underset{N\to\infty}{\longrightarrow}}

\newcommand{\E}{\op{\mathbb{E}}}

\newcommand{\R}{\mathbb{R}}
\newcommand{\C}{\mathbb{C}}
\newcommand{\n}{\mathbb{N}}

\newcommand{\ud}{\mathrm{d}}

\newcommand{\pro}{probability }

\newcommand{\f}{\frac}
\newcommand{\ff}{\frac{1}}
\newcommand{\lf}{\left}
\newcommand{\ri}{\right}

\newcommand{\st}{such that }

\newcommand{\lam}{\lambda}

\newcommand{\ti}{\times}

\newcommand{\vfi}{\varphi}
\newcommand{\ste}{\, ;\, }
\newcommand{\mc}{\mathcal }

\newcommand{\eps}{\varepsilon}

\newcommand{\bck}{\backslash}
\newcommand{\al}{\alpha}

\newcommand{\ovl}{\overline}

\newcommand{\bbm}{\begin{bmatrix}}
\newcommand{\ebm}{\end{bmatrix}}
\newcommand{\bes}{\begin{equation*}}
\newcommand{\ees}{\end{equation*}}
\newcommand{\be}{\begin{equation}}
\newcommand{\ee}{\end{equation}}
\newcommand{\beqy}{\begin{eqnarray}}
\newcommand{\eeqy}{\end{eqnarray}}
\newcommand{\beq}{\begin{eqnarray*}}
\newcommand{\eeq}{\end{eqnarray*}}
\newcommand{\one}{\mathbbm{1}}
\newcommand{\lto}{\longrightarrow}

\newcommand{\ie}{\emph{i.e. }}

\newcommand{\bpm}{\begin{pmatrix}}
\newcommand{\epm}{\end{pmatrix}}

\newcommand{\wt}{\widetilde}

\newcommand{\bpr}{\beg{proof}}
\newcommand{\epr}{\en{proof}}

\newcommand{\del}{\delta}

\newcommand{\pa}{\partial}

\newcommand{\ba}{\mathbf{a}}

\newcommand{\tvfi}{\tilde{\vfi}}
\newcommand{\tF}{\tilde{F}}

\newcommand{\sgn}{\mathrm{sgn}}

\newcommand{\ii}{\mathrm{i}}
\newcommand{\Gsc}{G_{\op{sc}}}

\theoremstyle{definition}

\long\def\symbolfootnote[#1]#2{\begingroup
\def\thefootnote{\fnsymbol{footnote}}\footnote[#1]{#2}\endgroup}

\allowdisplaybreaks

\title{Fluctuations of linear statistics of half-heavy-tailed random matrices}
\author{Florent Benaych-Georges and Anna Maltsev}



\thanks{FBG: MAP5,
Universit\'e Paris Descartes,
45, rue des Saints-P\`eres
75270 Paris Cedex 06, France. florent.benaych-georges@parisdescartes.fr. AM: Howard House, Departments of Mathematics, University of Bristol,
Bristol BS81SN, UK.
annavmaltsev@gmail.com. A.M. acknowledges the support of the Leverhulme Trust Early Career
Fellowship (ECF 2013-613).}

\begin{document}
\maketitle
\beg{abstract}
In this paper, we consider a Wigner matrix $A$ with entries whose cumulative distribution decays as $x^{-\alpha}$ with $2<\alpha<4$ for large $x$. We  are interested in the fluctuations of the linear statistics $N^{-1}\operatorname{Tr} \varphi(A)$, for some nice test functions $\varphi$. The behavior of such fluctuations has been understood for both heavy-tailed matrices (i.e. $\alpha < 2$) in \cite{ACFTCL} and light-tailed matrices (i.e. $\alpha > 4$) in \cite{bai-silver-book}. This paper fills in the gap of understanding it for $2 < \alpha < 4$. We find that while linear spectral statistics for heavy-tailed matrices have fluctuations of order $N^{-1/2}$ and those for light-tailed matrices have fluctuations of order $N^{-1}$, the linear spectral statistics for half-heavy-tailed matrices exhibit an intermediate $\alpha$-dependent order of $N^{-\alpha/4}$.
 \en{abstract}

%arxiv abstract:
%We consider a Wigner matrix $A$ with entries tail decaying as $x^{-\alpha}$ with $2<\alpha<4$ for large $x$ and study  fluctuations of linear statistics $N^{-1}\operatorname{Tr}\varphi(A)$. The behavior of such fluctuations has been understood for both heavy-tailed matrices (i.e. $\alpha < 2$) and light-tailed matrices (i.e. $\alpha > 4$). This paper fills in the gap of understanding for $2<\alpha<4$. We find that while linear spectral statistics for heavy-tailed matrices have fluctuations of order $N^{-1/2}$ and those for light-tailed matrices have fluctuations of order $N^{-1}$, the linear spectral statistics for half-heavy-tailed matrices exhibit an intermediate $\alpha$-dependent order of $N^{-\alpha/4}$.

\section{Introduction}
Let $A=[a_{ij}]$   be an $N\ti N$ Hermitian random matrix whose entries are i.i.d. and let $\lam_1, \ld, \lam_N$ be its eigenvalues. It is well known that if the entries of $A$ are duly renormalized, then for any continuous bounded test function $\vfi$, the random variable \be\la{2151417h}\ff{N}\Tr \vfi(A)=\ff{N}\sum_{i=1}^N \vfi(\lam_i)\ee has a deterministic limit, which is equal to the integral of $f$ with respect to the limit spectral distribution of $A$, namely the semicircle law when the entries have at least a second moment \cite{bai-silver-book,agz} and
different distributions depending on $\alpha$ if the entries are heavy-tailed with exponent $\al\in (0,2)$ (see \cite{BAGheavytails,CB,BCC}). The rate of convergence of the random variables of \eqre{2151417h} to its limit is not usually $\ff{\sqrt{N}}$, as i.i.d. $\lam_i$'s would give. In particular, if the entries of $A$ have a fourth moment, then the fluctuations of  $\ff{N}\Tr \vfi(A)$
around its expectation have order $\ff{N}$ (see \cite{bai-silver-book,KKP96,sinai,johansson, BaiYaoBernoulli2005, BAI2009EJP, lytova, MShcherbina11}). On the other hand, if the entries are heavy-tailed with exponent $\al\in (0,2)$ or Bernoulli with parameter of order $N^{-1}$, then  the fluctuations of  $\ff{N}\Tr \vfi(A)$
around its expectation have order $N^{-1/2}$ \cite{ACFTCL}. This difference of order in the fluctuations is due to the fact   that when the entries of $A$ have enough moments,  the eigenvalues of $A$ fluctuate very little,
as studied by Erd\"os, Schlein, Yau, Tao, Vu and their co-authors, who analyzed their rigidity in e.g. \cite{ESY2,EYY,Tao-Vu_0906.0510}. On the other hand, the heavier the tails the more similar to a sparse matrix  the (renormalized) matrix $A$ is, and the more independently its eigenvalues behave.

A finite fourth moment means that for large $x$, $\p(|a_{ij}|>x)\approx x^{-\al}$   with $\al>4$, whereas  heavy-tailed entries with exponent $\al\in (0,2)$ correspond precisely to $\p(|a_{ij}|>x)\approx x^{-\al}$ with  $\al\in (0,2)$. In this text, we
  fill in the gap of understanding the role of $\al$ in the fluctuations linear spectral statistics: when $\al\in (2,4)$, we
prove a central limit theorem for $\ff{N}\Tr \vfi(A)$ in the case where $\vfi$ is a sum of resolvent functions, it appears that the order of the fluctuations, in this case, is $N^{-\al/4}$. This completes the picture, summarized in Table \re{label_tableau}.

  \begin{table}[ht]
\begin{tabular}{|c|c|c|c|}
\hline
 & \; $\al <2$ \; & \; $2<\al<4$ \; & \; $\al>4$ \; \\
 \hline
\hline
Order of the fluctuations of \eqre{2151417h} &  $N^{-1/2}$ & $N^{-\al/4}$ & $N^{-1}$ \\
\hline
\end{tabular}
\caption{Orders of the fluctuations of the r.v. of \eqre{2151417h} around its expectation  as a function of the exponent $\al$ \st $\p(|a_{ij}|>x)\approx x^{-\al}$ for $x$ large.}\la{label_tableau}
\end{table}
Viewed in the light of concentration inequalities for linear spectral functionals  of random matrices, random matrices with half-heavy tailed entries interpolate between two extreme regimes, as shown in Table \re{label_tableau}
:\bgt
\ite Using only the independence of the entries,   (see  \cite[Lem. C.1]{alice-charles-HT} or   \cite{PasturBook,BCC,BCC2})     for any bounded function $\vfi : \R\to \R$ with finite total variation,  we have, for any $\del>0$,  \be\la{conc1}\p\lf( \lf|\Tr \vfi(A)-\E\Tr \vfi(A)\ri|\ge \del\ri) \ \le \ Ce^{-c\f{\del^2}{N\|\vfi\|^2_{\op{TV}}}},\ee  which proves that $$\sqrt{N}\lf( \ff{N}\Tr \vfi(A)-\E[\ff{N}\Tr \vfi(A)]\ri)$$ is bounded in \pro and explains why    the order of  the fluctuations of \eqre{2151417h} cannot be larger than $N^{-1/2}$,
\ite  In the case where the entries of  $A$   are independent and  satisfy a Log-Sobolev inequality (for example in the GO(U)E case), then
for any Lipschitz function $\vfi : \R\to \R$, we have, by \cite[Th. 2.3.5]{agz} (see also \cite{GZ2000,ledoux-amsbook, BLM}), \be\la{conc2}\p\lf( \lf|\Tr \vfi(A)-\E\Tr \vfi(A)\ri|\ge \del\ri) \ \le \ Ce^{-c\f{\del^2}{\|\vfi'\|^2_{\infty}}},\ee which proves that $$N\lf( \ff{N}\Tr \vfi(A)-\E[\ff{N}\Tr \vfi(A)]\ri)$$ is bounded in probability. It explains why the order of  the fluctuations of functionals as \eqre{2151417h} cannot be larger  than $N^{-1}$ in the case of matrices with independent Log-Sobolev entries.
\ent
Equation \eqre{conc1} shows that the case $\al<2$ corresponds to the largest possible fluctuations order in  \eqre{2151417h}. On the other hand,  \eqre{conc2} proves that in the
Gaussian case (this has been extended by Bai \emph{et. al.} to the case  $\al \ge  4$), the actual order is $N^{-1}$ (to be more precise, \eqre{conc2} only gives an upper-bound for this order, but one can easily check, using, for example, $\vfi(\lam)=\lam$ or $\vfi(\lam)=\lam^2$, that $N^{-1}$ is actually  the right order).  The case $2<\al<4$ is an intermediate case, where  concentration inequalities neither allow to guess the order of the fluctuations, nor allow to extend fluctuation results from a first class of test functions to a wider classer (as was done for example in  \cite{ACFTCL}).

\noindent{\bf Notation.} Here, $A\ll B$ means  that $A/B\lto 0$ as $N\to\infty$,  $A\sim B$ means  that $A/ B\lto 1$ as $N\to\infty$ (or $x\to\infty$ if $x$ is the underlying variable) and $A=O(B)$ means that $A/B$ stays bounded as $N\to\infty$. If $A$ and $B$ are random variables, similar notations are used in probability. Also, for $\al\in \R$, the function $z\mapsto z^\al$ is defined on $\C\bck\R^-$ thanks to the determination of the argument which is null on $(0, +\infty)$. At last, for functions of a complex variable $z=x+\ii y$, we use the classical notation $\pa_z=(\pa_x-\ii\pa_y)/2$.

\section{Main result}
%We say that a random variable $x_{ij}$ is half-heavy-tailed    with parameter $\alpha$ if \st for a certain $\al\in (2,4)$ and a certain $c >0$, \be\la{tailX} \p(|x_{ij}|>x)\sim \f{c}{\Ga(\al+1)}x^{-\al}\ee for large $x$.
Let us consider a random real symmetric or Hermitian  matrix $$\ds A=[a_{ij}]_{1\le i,j\le N}=\left[\f{\ds x_{ij}}{\sqrt{N}}\right]_{1\le i,j\le N},$$
where one of two conditions holds: either
\begin{enumerate}
\item[(a)]\label{real}(real case) $x_{ij}$'s, $1\le i\le j$, are i.i.d. real random variables with mean $0$ and variance $1$ \st for a certain $\al\in (2,4)$ and a certain $c >0$, as $x\to+\infty$, \be\la{tailX} \p(|x_{ij}|>x)\sim \f{c}{-\Ga(1-\al/2)}x^{-\al},
\ee
or
\item[(b)]\label{complex}(complex case) $x_{ij} = x_{ij}^R/\sqrt 2 + \ii x_{ij}^I/\sqrt 2$ for $1 < i < j$ and $x_{ii} = x_{ii}^R$ where $x_{ij}^I$ and $x_{ij}^R$ are i.i.d. real symmetric random variables with mean $0$ and variance $1$ that satisfy \eqre{tailX}.
\end{enumerate}

Our main theorem is the following.  % Let $\mc{E}$ denote the set of functions $\R\to\C$ which are linear combinations of functions of the type $\vfi_z:\lam\longmapsto \ff{z-\lam}$, with $z\in \C\bck\R$. Note that by the Stone-Weierstrass theorem,  $\mc{E}$ is dense in the set of continuous functions $\R\to\C$ with null limit at infinity endowed by the norm $\|\vfi\|_\infty:=\sup_{\lam\in \R} |\vfi(\lam)|$.

\beg{Th}\la{maintheorem}
For $$G(z):=(z-A)^{-1} $$  with $A$ as above the process $$\lf(\ff{N^{1-\al/4}} (\Tr G(z)-\E\Tr G(z))\ri)_{z\in \C\bck\R}$$ converges to a complex Gaussian centered process $(X_z)_{z\in \C\bck\R}$ with covariance  defined by the fact that $X_{\ovl{z}}=\ovl{X_z}$ and that for any $z,z'\in \C\bck\R$, $\E[X_zX_{z'}]=C(z,z')$, for \beq C(z,z')&:=&\iint_{t,t'>0} \pa_z\pa_{z'}\bigg\{[(K(z,t)+K(z',t'))^{\al/2}-(K(z,t)^{\al/2}+K(z',t')^{\al/2})]\\ &&\exp\lf(\sgn_z \ii tz-K(z,t)+\sgn_{z'} \ii t'z'-K(z',t')\ri)\bigg\}\f{c\,\ud t\,\ud t'}{2tt'}\eeq
where $c$ and $\al$ are as in \eqre{tailX},  $\sgn_z:=\sgn(\Im z)$ and $K(z,t):=\sgn_z \ii t\Gsc(z)$, $\Gsc(z)$  being the Stieltjes transform of the semicircle law with support $[-2,2]$.
\en{Th}

\beg{rmk}This theorem proves Gaussian convergence for any random variable of the form  $$\ff{N^{1-\al/4}} (\Tr \vfi(A)-\E\Tr \vfi(A)),$$ where $\vfi$ is a function of the type $$\qquad \qquad \qquad \qquad \vfi(\lam)=\sum_{j=1}^p \f{c_j}{z_j-\lam},\qquad \qquad (p\ge 1,\; c_1, \ld, c_p\in \C,\; z_1, \ld, z_p\in \C\bck\R).$$ The functions $\vfi$ of this type span (by closure) some larger sets of  functions (by, e.g.
the Stone-Weierstrass theorem,  the Cauchy formula, or the Helffer-Sj\"ostrand formula).   However, the lack of error control in approximating $\vfi$ (due to the fact that in our case the Log-Sobolev concentration inequality \eqre{conc2} is not true and the general concentration inequality \eqre{conc1} is not sharp enough) prevents us from extending our theorem to a larger class of test functions, as was done from resolvent functions to wider classes in e.g. \cite{bai-silver-book,ACFTCL,MShcherbina11}. As far as applying    \cite[Propo. 1]{MShcherbina11} is concerned, the problems come first from the fact that we truncate the entries of the matrices to upper-bound the variance of $\Tr (E+\ii\eta -A)^{-1}$ and second from the fact that this variance does not decay enough as $|E|$ grows.
\en{rmk}

The remainder of the paper consists of the proof of Theorem \ref{maintheorem}. In Section \ref{truncation}, we truncate the random variables appropriately, and centralize in the real case  (in the complex  case, centralization is automatic due to our assumption of symmetry). In Section \ref{martingale}, we restate our problem in terms of a martingale approach and cite relevant martingale convergence theorem. In Section \ref{removing}, we show that off-diagonal terms of the resolvent can be neglected in further calculations. Lastly, in Section \ref{calculation} we show that the diagonal terms of the resolvent yield the desired formula for the covariance, using a lemma proved in Section \ref{concentration} that allows us to approximate the diagonal elements of the resolvent by the Stietjes transform of the spectral measure.\\ \\ \\

\section{Proof of Theorem \re{maintheorem}}\label{s:symmetric}
\subsection{Truncation and recentralization of the entries} \label{truncation}

Let $A,B$ be any $N\ti N$ matrices, $z\in \C\bck\R$, and  \be G_B(z):=\ff{z-B} .\ee
Then \be
\la{f}G_B(z)=\ff{z-A}+\ff{z-A}(B-A)\ff{z-B}\ee and we have that
$$|\Tr( G_B(z)-G_A(z))|\le 2|\Im z|^{-1}\op{rank}(B-A).$$
Thus for fixed $z$,
$$\op{rank}(B-A)\ll N^{1-\al/4}\implies |\Tr( G_B(z)-G_A(z))|\ll N^{1-\al/4}.$$
Let us consider the case of real symmetric matrices and let $\mu_N = \E x_{ij}\one_{|x_{ij}|\le N^\beta}$ (for an exponent $\beta$ which will be specified later). First we estimate the absolute value of $\mu_N$. As $x_{ij}$ is centered, \be\la{312151}\mu_N:=\E x_{ij}\one_{|x_{ij}|> N^\beta},\ee so for $N$ large enough,
 \begin{equation}\begin{split}
0 \le |\mu_N| & \le  \E |x_{ij}| \one_{|x_{ij}| \ge N^{\beta}} = \int_{N^{\beta}}^{\infty}\p (|x_{ij}| > x)\ud x
 \, \le \, \f{2c}{-\Ga(\al/2-1)( \al-1)}(N^{\beta})^{-\al+1}
\end{split}\end{equation}
Let $B=[a_{ij}\one_{|x_{ij}|\le N^\beta} - \mu_N/\sqrt N]$. Subtracting $\mu_N/\sqrt N$ from each matrix entry is a rank 1 perturbation. Then, as
$$\p(|x_{ij}|> N^\beta)\le  CN^{-\al\beta},$$
we have
$$\op{rank}(B-A) \le 1+ 2\sum_{i=1}^NX_i$$ where the $X_i$'s are independent Bernoulli r.v. with parameters $$\p(X_i=1)\;=\;1-(1-CN^{-\al\beta})^i$$
and 1 is added for the rank 1 perturbation of shifting each entry by $\mu_N/\sqrt N$.
 In order to upper-bound $\op{rank}(B-A)$ with high \pro thanks to Bennett's inequality or Lemma 5.7 in \cite{floSandrine}, we compute the sum of these parameters:\beq \sum_{i=1}^N1-(1-CN^{-\al\beta})^i&=& N -( 1-CN^{ -\al \beta }) \f{ 1- ( 1-CN^{ -\al \beta} )^N} {CN^{ -\al \beta}}.
\eeq
If $\al\beta>1$,  \beq (1-CN^{-\al\beta})^N&=&\exp N\log(1-CN^{-\al\beta})\\
&=&1-CN^{1-\al\beta}+\f{C^2}{2}N^{2(1-\al\beta)}+o(N^{2(1-\al\beta)})\eeq
so \beq \sum_{i=1}^N1-(1-CN^{-\al\beta})^i&=&N-(1-CN^{-\al\beta})\f{CN^{1-\al\beta}-\f{C^2}{2}N^{2(1-\al\beta)}+o(N^{2(1-\al\beta)})}{CN^{-\al\beta}}\\
&=&N-(1-CN^{-\al\beta})(N-\f{C}{2}N^{2-\al\beta}+o(N^{2-\al\beta})) \\
&=&\f{C}{2}N^{2-\al\beta}+o(N^{2-\al\beta})
\eeq
Thus, by Bennett's inequality or Lemma 5.7 in \cite{floSandrine},  as soon as $\al\beta>1$, we know that $\op{rank}(B-A)$ has order at most $N^{(2-\al\beta)_+}$ (\ie for any $\eps>0$, $N^{-(2-\al\beta)_+-\eps}\op{rank}(B-A)$ tends in \pro to zero).
So one can replace $A$ by $B$ as long as $$(2-\al\beta)_+<1-\al/4,$$\ie   $$\beta>\f{2-({1-\al/4})}{\al}= \frac{1}{\alpha}+ \frac{1}{4}.$$
Furthermore we want to renormalize our new truncated centered random variables to have variance 1, so we let
\be\la{312152}\sigma_N := \sqrt{\E(x_{ij} \one_{|x_{ij}| < N^{\beta}} - \mu_N)^2}
\ee
Noting that
\[
 \E(\one_{|x_{ij}| >N^{\beta}}x_{ij}^2) = \int_{N^{2\beta}}^{\infty} \p(|x_{ij}|^2 > x) \ud x\, \sim \,  \frac{c}{\Gamma(\alpha + 1)(1-\alpha/2)N^{\beta(\alpha-2)}}  , \]
we see that
\[ \si_N^2-1=
  O(N^{\beta(2 - \alpha)}).
\]
and thus $\si_N = 1 + O(N^{\beta(2 - \alpha)})$. We keep the variance of the entries as $\sigma_N^2/N$ throughout.

In the complex  case, subtracting the mean from each matrix entry is no longer a rank 1 perturbation, so this argument will no longer work. This is the reason why, in the complex case, we only consider random variables which are symmetric so that we can truncate and still retain a $0$ mean.
\beg{lem}\la{215141}Let  $\epsilon>0$,  $\beta = \frac{1}{\alpha}+\frac{1}{4}+\epsilon$  and $$\tilde{a}_{ij}:=(x_{ij}\one_{|x_{ij}| < N^{\beta}}- \mu_N)/\sqrt N$$ in the real case   and $$\tilde{a}_{ij}:= (x_{ij}^R\one_{|x_{ij}^R| < N^{\beta}} + \ii x_{ij}^I\one_{|x_{ij}^I| < N^{\beta}} )/\sqrt N$$ in the complex case. Then there is a constant $C$ depending only on the distribution of the $x_{ij}$'s \st
 \bgt\ite[(i)] the $\tilde{a}_{i,j}$'s are i.i.d., centered, with variance $\sigma_N^2/N$,
 \ite[(ii)] $N^{3/2}\E [|\tilde{a}_{i,j}|^3]\le CN^{\left(\frac{1}{\alpha}+ \frac{1}{4} + \epsilon\right)(3-\al)_+} $,%\leq \max\{C, N^{\frac{3}{4} - \frac{\alpha}{16} - \frac{ \alpha^2} {16}+\epsilon(3-\al)}\}$
 \ite[(iii)] $N^2\E [|\tilde{a}_{i,j}|^4]\le CN^{\beta(4-\al)}=CN^{\frac{4}{\alpha}-\frac{\alpha}{4}+\epsilon(4-\al)}$,
 \ite[(iv)] for any $\lam\in \C$ \st $\Im \lam\le  0$, \be\la{2071216h23} \phi_N(\lam):=\E[e^{-\ii\lam|\tilde{a}_{ij}|^2}]=1-\f{\ii\lam\sigma_N^2}{N} + c\f{(\ii\lam \sigma_N^2)^\f{\al}{2}}{N^{\f{\al}{2}}}+\f{|\lam\sigma_N|^{\f{\al}{2}}}{ N^{\f{\al}{2}}}\eps_N(\ii\lam \sigma_N^2/N),\ee
where the function  $\eps_N(z)$ is analytic in $z$ on $\{\Re(z)>0\}$,  bounded uniformly in $(N,z)\in \n\ti \mc{K}$ for  any compact $\mc{K}\subset\{ \Re(z)\ge 0\}$ and $\lim_{z\to 0}\eps_N(z)=0$  uniformly in $N$. \ent\en{lem}
 \begin{proof}
 (i) is true by the definition of $\mu_N$ and $\si_N$ at \eqre{312151} and \eqre{312152}, (ii) and (iii) are easy computations, relying on the fact that  for $f:\R^+\to \R^+$ increasing  and $X$ a positive random variable, $\E[f(X)]=f(0)+\int_0^{+\infty} f'(t)\p(X\ge t)\ud t$.
 Lastly, (iv) follows from  \cite[Th. 8.1.6]{Bingham-Goldie-Teugels} for non-truncated entries. To centralize,  note that
  \[
 \p(|x_{ij} - \mu_N|\ge x) = \p(x_{ij} \ge  x + \mu_N) + \p(x_{ij}\le  -x + \mu_N \; \trm{and} \; x_{ij} \le  \mu_N). \]
 Recalling that $\mu_N\to  0$ we get from \eqref{tailX} that  for a certain constant $C>0$,
\[
\f{c}{C-\Ga(\al/2-1)}(x-\mu_N)^{-\al}
\le \p(|x_{ij} - \mu_N|>x) \le  \f{Cc}{-\Ga(\al/2-1)}(x+\mu_N)^{-\al}
\]
Since $\mu_N \rightarrow 0$, the random variable $x_{ij} - \mu_N$ will also satisfy (\ref{tailX}) for large $N$, and therefore (iv) holds for the shifted entry $x_{ij} - \mu_N$.
\end{proof}

From now on, we suppose that each $a_{ij}$ has been replaced by the $\tilde{a}_{ij}$ of the previous lemma,
 for \be\la{1751417h}\beta:=\frac{1}{\alpha}+ \frac{1}{4} +\epsilon,\ee for $\epsilon>0$ that can be chosen as small as needed.   By a slight abuse of notation,    we  still denote this random variable  by $a_{ij}$  and we henceforth assume the conclusions of Lemma \re{215141} to be true for the $a_{ij}$'s.

% \textbf{From now on, we suppose Lemma \re{215141} to hold.\\ \\}

\subsection{Martingale approach}\label{martingale}
%Let $A$ be a heavy-tailed random matrix, whose eigenvalues are denoted by $\lam_i$.
We want to prove that for certain test functions $\vfi$ (namely the linear combinations of functions of the type $\lam\in \R\longmapsto \ff{z-\lam}$, with $z\in \C\bck\R$), $$M(\vfi,N):=\ff{N^{1-\al/4}}(\Tr\vfi(A)-\E[\Tr\vfi(A)]) $$
converges in distribution to a certain Gaussian distribution.
We will use  Theorem \re{thconvmart} for $M(\vfi,N)$, with $M_N(N)=M(\vfi,N)$ and $\mc{F}_k(N):=\si(x_{i,j}\ste i\le k \trm{ and }j\le k).$

Then, denoting $\E[\,\cdot\,|\mc{F}_k]$ by $\E_k$, the random variable $Y_k(N)$ of Theorem \re{thconvmart} is  $$Y_k=Y_k(N)=\ff{N^{1-\al/4}}(\E_k-\E_{k-1})(\Tr \vfi(A)).$$

Let $A^{(k)}$ be the $N-1$ by $N-1$ matrix obtained by removing the $k$th row and column of $A$. Then $(\E_k-\E_{k-1})(\Tr \vfi(A^{(k)}))=0$, hence$$Y_k=\ff{N^{1-\al/4}}(\E_k-\E_{k-1})(\Tr \vfi(A)-\Tr \vfi(A^{(k)})).$$

%{\bf Thus it appears that using  Theorem \re{thconvmart} in our case comes down to sum up the dependencies of our statistics on the columns of $A$.}

Note first that by the interlacing property between the spectrums of $A$ and $A^{(k)}$, when $\vfi$ has finite total variation, we have \be\la{consinterlac}|\Tr \vfi(A)-\Tr \vfi(A^{(k)})|\le \|\vfi\|_{\op{TV}}.\ee As a consequence, $|Y_k|\le \f{ \|\vfi\|_{\op{TV}}}{N^{1-\al/4}}$ and the $L(\eps, N)$ of  Theorem \re{thconvmart} is null for $N$ large enough.

Hence it remains to prove that $$\sum_{k=1}^N \E_{k-1}(Y_k^2)$$ and  $$\sum_{k=1}^N \E_{k-1}(|Y_k|^2)$$ have finite deterministic limits that agree with the limit covariance of Theorem \re{maintheorem}.

  By linear combination, it suffices to prove that for  any $z,z'\in\C\bck\R$, for   $$Y_k:=\ff{N^{1-\al/4}}(\E_k-\E_{k-1})(\Tr \ff{z-A}-\Tr \ff{z-A^{(k)}}) $$ and $$Y_k':=\ff{N^{1-\al/4}}(\E_k-\E_{k-1})(\Tr \ff{z'-A}-\Tr \ff{z'-A^{(k)}}),$$ we have $$\sum_{k=1}^N \E_{k-1}( Y_kY_k' )$$ converges in \pro to a deterministic constant which agrees with the limit covariance of Theorem \re{maintheorem}.

Note that $$\sum_{k=1}^N \E_{k-1}( Y_kY_k' )= \ff{N}\sum_{k=1}^N  N\E_{k-1}( Y_kY_k' ) = \int_{u=0}^1
N\E_{\lceil Nu\rceil-1}( Y_{\lceil Nu\rceil}Y_{\lceil Nu\rceil}' )
\ud u,$$%\qquad \trm{ for $k:=\lceil Nu\rceil$},$$
hence we shall prove that for any $u\in (0,1)$, as $N\to\infty$ and $k\to\infty$ with $k/N\to u$, we have  \be\la{DefC(u)513}N\E_{k-1}( Y_kY_k' )\lto 2uC(z,z'),\ee with $C(z,z')$ the function defined in Theorem \re{maintheorem}.

Note also that for $G^{(k)}(z):=(z-A^{(k)})^{-1}$, by  \eqre{307121},
\be\la{30712121114} Y_k=\ff{N^{1-\al/4}}(\E_k-\E_{k-1})\f{1+\ba_k^*(G^{(k)}(z))^{2}\ba_k}{z-a_{kk}-\ba_k^*G^{(k)}(z)\ba_k}
 \ee
where $\ba_k$ is the $k$th column of $A$ without the diagonal term.

\subsection{Removing the off-diagonal terms}\label{removing}
\beg{propo}\la{85141}   Let us define $$\tilde{Y}_k:=
\ff{N^{1-\al/4}}(\E_k-\E_{k-1})\f{1+\ba_k^*(G^{(k)}(z))^2_{\op{diag}}\ba_k}{z-\ba_k^*G^{(k)}(z)_{\op{diag}}\ba_k}
$$ and $$\tilde{Y}'_k:=
\ff{N^{1-\al/4}}(\E_k-\E_{k-1})\f{1+\ba_k^*(G^{(k)}(z'))^2_{\op{diag}}\ba_k}{z'-\ba_k^*G^{(k)}(z')_{\op{diag}}\ba_k},
$$ where for a matrix $M$, $M_{\op{diag}}$ denotes the diagonal matrix obtained from $M$ by setting all its non-diagonal entries to zero.  Then \be\la{276132}\sum_{k=1}^N\E_{k-1}[ {Y}_k {Y}_k']-\E_{k-1}[\tilde{Y}_k\tilde{Y}_k']\ee converges in \pro to $0$.
\en{propo}

\bpr We define, for $z\in \C\bck \R$,   $$F_k:=\log|z-a_{kk}-\ba_k^*G^{(k)}(z)\ba_k|^2\qquad  F'_k:=\log|z'-a_{kk}-\ba_k^*G^{(k)}(z')\ba_k|^2$$
$$\tF_k:=\log|z- \ba_k^*G^{(k)}(z)_{\op{diag}}\ba_k|^2\qquad  \tF'_k:=\log|z'-\ba_k^*G^{(k)}(z')_{\op{diag}}\ba_k|^2.$$
These functions are well defined because for   any $z\in \C\bck\R$, $$\Im z \ti \Im(-\ba_k^*G^{(k)}(z)\ba_k)>0\qquad;\qquad\Im z \ti \Im(-\ba_k^*G^{(k)}(z)_{\op{diag}}\ba_k)>0$$which implies that \be\la{41220151} \beg{cases} \Im z>0\implies \Im(z-a_{kk}-\ba_k^*G^{(k)}(z)\ba_k), \Im(z-a_{kk}-\ba_k^*G^{(k)}(z)_{\op{diag}}\ba_k)>\Im z\\
\Im z<0\implies \Im(z-a_{kk}-\ba_k^*G^{(k)}(z)\ba_k), \Im(z-a_{kk}-\ba_k^*G^{(k)}(z)_{\op{diag}}\ba_k)<\Im z\end{cases}
 ,\ee so that the argument of the $\log$ cannot vanish.

% For many further uses, we state for once the following lemma:
% \beg{lem}Let $(X(z,x))_{x\in \C\bck$ be
% \en{lem}

Using the fact that, for an analytic function $f$ defined on $\C\bck\R$ and taking values in $\C\bck\R$ \st $f(\ovl{z})=\ovl{f(z)}$, $$\pa_z \log|f(z)|^2=\pa_z \log(f(z)f(\ovl{z}))=\pa_z\log(f(z))=\f{f'(z)}{f(z)},$$we have
$$ Y_k= \ff{N^{1-\al/4}}(\E_k-\E_{k-1})\pa_z F_k
= \pa_z\ff{N^{1-\al/4}}(\E_k-\E_{k-1})F_k,
$$where in the second equality, we used \eqre{41220151} to   notice the   $F_k$ is  bounded uniformly in the randomness and in $z$ as $z$ varies in any compact subset of $\C\bck\R$. Of course analogous formulas hold for $Y_k'$, $\tilde{Y}_k$ and $\tilde{Y}_k'$.
 Thus, commuting conditional expectation and derivative again for the same reason, we have
\beq&&\sum_{k=1}^N\E_{k-1}[ {Y}_k {Y}_k']-\E_{k-1}[\tilde{Y}_k\tilde{Y}_k'] =\\ &&\pa_z\pa_{z'}\sum_{k=1}^NN^{-2+\al/2}\{\underbrace{\E_{k-1}[ (\E_k-\E_{k-1})F_k (\E_k-\E_{k-1})F_k']}_{\ds :=\vfi_k}-\underbrace{\E_{k-1}[(\E_k-\E_{k-1})\tF_k (\E_k-\E_{k-1})\tF_k']}_{\ds :=\tvfi_k}\}\eeq
hence by the Cauchy inequalities for holomorphic functions, it suffices to prove that
%there is a finite constant $C$ independent on $N$ and $k$ and
 uniformly on $k,z,z'$ (as $z,z'$ stay at a macroscopic  distance   from the real line) we have   \be\la{1851423h09}|\vfi_k-\tvfi_k|\ll  N^{-1}N^{2 - \alpha/2}.\ee

Let
\begin{equation}\label{d:etak}\eta_k:=a_{kk} + \ba_k^*G^{(k)}(z) \ba_k - \ba_k^*G^{(k)}(z)_{\op{diag}} \ba_k = a_{kk} + \sum_{i\ne j} G^{(k)}(z)_{ij} \ovl{\ba_k(i)} \ba_k(j),\end{equation} and   \begin{equation} \eps_k\ :=\ F_k-\tF_k\  = \log |1-\eta_k(z- \ba_k^*G^{(k)}(z)_{\op{diag}}\ba_k)^{-1}|^2.\end{equation}
We also define $\eta_k'$ and $\eps_k'$ in the same way with $z$ instead of $z'$.

By Lemma \re{1912131}, we have $$\vfi_k=\E_{k-1}[\E_k(F_k)\E_k(F_k')]-\E_{k-1}(F_k)\E_{k-1}(F_k')$$ and the analoguous equality holds for $\tvfi_k$. Thus using the formulas $F_k=\tF_k+\eps_k$ and $F_k'=\tF'_k+\eps'_k$, we have
\beqy\la{185141} \vfi_k-\tilde{\vfi}_k&=& \E_{k-1}[\E_k(\tilde{F}_k)\E_k(\eps_k')]+\E_{k-1}[\E_k(\eps_k)\E_k(\tilde{F}'_k)]+\E_{k-1}[\E_k( \eps_k)\E_k( \eps_k')]\\ \nonumber&&-\E_{k-1}(\tilde{F}_k)\E_{k-1}(\eps_k')-\E_{k-1}(\eps_k)\E_{k-1}(\tilde{F}_k)-\E_{k-1}(\eps_k)\E_{k-1}(\eps_{k}')
\eeqy
%hence\be\la{276131}\| \vfi-\tilde{\vfi} \|_{L^1}\le 2\|\tilde{X}\|_{L^\infty}\|\eps_Y\|_{L^1}+2\|\tilde{Y}\|_{L^\infty}\|\eps_X\|_{L^1}+2\|\eps_X\|_{L^2}\|\eps_Y\|_{L^2}.\ee

Let $\E_{\ba_k}$ denote the expectation with respect to the randomness of the $k$th row $\ba_k$ of $A$. We have  \be\la{17414}\E_{k-1}(\,\cdot\,)=\E_k\E_{\ba_k}(\,\cdot\,)=\E_{\ba_k}\E_k(\,\cdot\,),\ee so that \eqre{185141} can be rewritten as
\beqy\la{185142} \vfi_k-\tilde{\vfi}_k&=& \E_{\ba_k}[\E_k(\tilde{F}_k)\E_k(\eps_k')]+\E_{\ba_k}[\E_k(\eps_k)\E_k(\tilde{F}'_k)]\\ \nonumber\\\nonumber&&+\E_{\ba_k}[\E_k( \eps_k)\E_k( \eps_k')]-\E_k\E_{\ba_k}(\tilde{F}_k)\E_k\E_{\ba_k}(\eps_k')\\ \nonumber\\\nonumber&&-\E_k\E_{\ba_k}(\eps_k)\E_k\E_{\ba_k}(\tilde{F}_k)-\E_k\E_{\ba_k}(\eps_k)\E_k\E_{\ba_k}(\eps_{k}')\\\nonumber
\\\nonumber&=:&  T_1 + T_2 + T_3 + T_4 + T_5 + T_6
\eeqy
From now on, $C$ will denote   a finite constant  (that will change from line to line) depending uniformly in $k,z,z'$ as $z,z'$ stay at any positive distance away from the real line.
\beg{lem}\label{l:eta24}We have
%\ite[(i)]  for $P$  a product with  an odd number of factors chosen among  $\eta_k$, $\ovl{\eta_k}$, $\eta_k'$, $\ovl{\eta_k'}$ is centered with respect to the partial expectation $\E_{\ba_k}$ and for $f$ a bounded function,$$\E_{\ba_k}[P\ti f(\ba_k^*G^{(k)}(z)_{\op{diag}}\ba_k,\ba_k^*G^{(k)}(z')_{\op{diag}}\ba_k]=0,$$
 $\E_{\ba_k}[|\eta_k|^2]\le CN^{-1},$
for $\epsilon$ as in \eqre{1751417h}.
Besides, the same bounds hold if one replaces $\eta_k$ by  $\eta_k'$.
\end{lem}

\bpr %(i) is due to the fact that the $a_{ij}$'s are symmetrically distributed.
Let us denote $G=G^{(k)}(z)$ and write $\E$ instead of $\E_{\ba_k}$ for short. Note that \be\la{4121518h55}\Tr GG^*\le CN\qquad;\qquad \|G\|\le C.\ee Inequality  (i) follows from  Lemma \re{2151411h50}, which allows to claim that   $$
	\E_{\ba_k} [| \eta_k|^2] \le  CN^{-1}.
$$	\epr

 %By \eqre{185141},
 To conclude the proof  of  Proposition \re{85141}, we need to prove the upper bound \eqre{1851423h09}.
We will use the expression of $\vfi_k-\tvfi_k$ given at  \eqre{185142}, denoting   the six terms of its RHS by $T_1, \ld, T_6$. We will show that $|T_1|, |T_2|, |T_3|, |T_4|, |T_5|,$ and $ |T_6| \ll N^{1 - \alpha/2}$.

%Here we use the assumption that the entries are symmetrically distributed to avoid the cumbersome analysis of some of the lower order terms. We show how to handle them in Section \ref{s:nosymmetry}.

To make subsequent calculations less cumbersome to write, we introduce the notation $$J_k = \frac{1}{z-a_{kk}-\ba_k^*G^{(k)}(z)\ba_k}$$ and
$$J_{k, \diag} = \frac{1}{z - \ba_k^*G^{(k)}_{\diag}(z)\ba_k},$$
and correspondingly $J_k', J_{k, \diag}'$ with $z'$ instead of $z$.
Let furthermore
\[
J_{k, \Tr} := \frac{1}{-z - \frac{1}{N} \Tr G^{(k)}(z)}
\]
and $E$ be given by
\[
E := \ba_k^* G^{(k)}(z)_{\op{diag}} \ba - \frac{1}{N}\Tr G^{(k)}(z) = \sum_j G^{(k)}_{jj} \left(|\ba_k(j)|^2 - \frac{1}{N}\right)
\]
so that
\[
J_{k, \diag} =J_{k, \Tr} + E J_{k, \diag} J_{k, \Tr}
\]
%First, note that if $a_{ij}$'s are symmetrically distributed, we have \be\la{1851414h51}\E_{\ba_k}[\eta_k \ti \trm{(any function of $J_{k, \diag}$ and $J_{k, \diag}'$)}]=0\ee and that the same holds for $\eta_k'$.

To find a bound on $|\E_{\ba_k}[\eps_k]|$, we bound $\E_{\ba_k}[\eps_k]$ and $\E_{\ba_k}[-\eps_k]$ separately. We notice that $|1-\eta_k J_{k, \diag}|$ and $|1 + \eta_k J_k|$ are reciprocals. Using Jensen's inequality, we find a bound on the former:
\begin{equation*}\begin{split}
\E_{\ba_k}[\eps_k] &= \E_{\ba_k}\log |1-\eta_k J_{k, \diag}|^2 \\
& \leq \log \E_{\ba_k}|1-\eta_k J_{k, \diag}|^2 \\
&
= \log \E_{\ba_k}(1 - 2 \Re (\eta_k J_{k, \diag}) + |\eta_k J_{k, \diag}|^2) \\
&
= \log \E_{\ba_k}(1 - 2 \Re (\eta_k J_{k, \Tr} + \eta_k E J_{k, \diag} J_{k, \Tr}
) + |\eta_k J_{k, \diag}|^2) \\
&
 \leq \log (1 + O(\E_{\ba_k}|\eta_k| \, |E|) + O(\E_{\ba_k}|\eta_k|^2))
\end{split}\end{equation*}
Here in the last line we have used that $J_{k, \Tr}$ is independent of $\ba_k$ and that $\E_{\ba_k} \eta_k = 0$ which gives that $\E_{\ba_k}\, \Re\, ( \eta_k J_{k, \Tr}) = 0$. We have also used that $J_k$, $J_{k, \diag}$, and $\ds J_{k, \Tr}$ are uniformly bounded to claim that $\E_{\ba_k}\eta_k E J_{k, \diag} J_{k, \Tr} \leq O(\E_{\ba_k}|\eta_k| \, |E| )$ and $\E_{\ba_k}|\eta_k J_{k, \diag}|^2 \leq O(\E_{\ba_k}|\eta_k|^2)$. Similarly, we show the same bound for $\E_{\ba_k}[-\eps_k]$:
\begin{equation*}\begin{split}
\E_{\ba_k}[-\eps_k] &= \E_{\ba_k}\log |1 + \eta_k J_k|^2 \\
& \leq \log \E_{\ba_k}|1 + \eta_k J_{k, \Tr} + \eta_k(\eta_k + E) J_{k, \diag} J_k|^2 \\
&
= \log (1 + 2 \E_{\ba_k} \Re (\eta_k J_{k, \Tr}) + O (\E_{\ba_k} |\eta_k| \, |E|) + O(\E_{\ba_k}|\eta_k|^2 )) \\
&
 \leq \log (1 + O (\E_{\ba_k} |\eta_k| \, |E|) + O(\E_{\ba_k}|\eta_k|^2 ))
\end{split}\end{equation*}
 By Cauchy-Schwarz and Lemma \ref{l:eta24} we get that
\beq
\E_{\ba_k}[|\eta_k| \, |E|] &\le&  \sqrt{\E_{\ba_k}|(\eta_k )|^2 \E_{\ba_k}| E|^2}
\\ &
 \le&  C N^{-\frac{1}{2}} N^{\frac{1}{2}\left(\frac{4}{\alpha} - \frac{\alpha}{4} - 1+\epsilon(4-\al)\right)}
\eeq
We note that for $2<\alpha<4$, it can be ensured that $\epsilon$ is small enough to obtain the desired inequality:
\beq
\E_{\ba_k}[|\eta_k| \, |E|] \ll N^{1 - \alpha/2}
\eeq
yielding the bound
\[
|\E_{\ba_k}[\eps_k]| \leq O (\E_{\ba_k} |\eta_k| \, |E|) \ll N^{1 - \alpha/2}
\]
since $O(\E_{\ba_k} |\eta_k| \, |E|) \ll1$ and therefore the logarithm is given by its Taylor series.
Hence  $|T_4|, |T_5|$ and  $|T_6| \ll CN^{1 - \alpha/2}$.

Next, after applying Cauchy-Schwartz to $T_3$, we look for a bound on $\E_{\ba_k}[|\E_k(\eps_k)|^2]$.
We note first that by Jensen's inequality
\[
\E_{\ba_k}[|\E_k(\eps_k)|^2] \leq \E_k \E_{\ba_k}[|(\eps_k)|^2]
\]
Let us partition the space of matrices as follows. We define the events
\be\la{2691310h}S := \{A: |1 + \eta_k J_k|^2 \geq 1 \}\qquad ;\qquad \tilde{S} := \{A: |1 - \eta_k J_{k, \diag}|^2 \geq 1 \} \ee
%and
%\[S^c := \{A: |1 + \eta_k J_k|^2 \geq 1 \}->S, S\to\tilde{S}.\]
and note that since $|1 - \eta_k J_{k, \diag}|$ and $|1 + \eta_k J_k|$ are reciprocals, we have that
\[
|\log |1 - \eta_k J_{k, \diag}|^2| = \one_{\tilde{S}} \log |1 - \eta_k J_{k, \diag}|^2 + \one_{S} \log |1 + \eta_k J_k|^2.
\]
Recall that for $x>0$, $\log (1 + x) < x$, so that
on $\tilde{S}$
\[
0 \leq (\log |1 - \eta_k J_{k, \diag}|^2) \one_{\tilde{S}} \leq (- 2\Re (\eta_k J_{k, \diag}) + |\eta_k J_{k, \diag}|^2)\one_{\tilde{S}}.
\]
Using the definition of $S$ and similar reasoning we get a similar bound on $\one_{S} \log |1 + \eta_k J_k|^2$ with $J_k$ instead of $-J_{k, \diag}$, which yields
\begin{equation}\label{uboundlog}
|\log |1 - \eta_k J_{k, \diag}|^2| \leq - \one_{\tilde{S}} (2\Re (\eta_k J_{k, \diag}) + |\eta_k J_{k, \diag}|^2) + \one_{{S}} ( 2\Re (\eta_k  J_k) + |\eta_k J_k|^2)
\end{equation}
Taking the expectation of $|\eps_k|^2$, and using that $|J_k|$ and $|J_{k, \diag}|$ are absolutely bounded we get that
\[
\E_k(\log |1 - \eta_k J_{k, \diag}|^2)^2 \leq C \E_k |\eta_k |^2 < CN^{-1}
\]
 where the last inequality follows by the Lemma \ref{l:eta24}. Here the $\epsilon$ of \eqre{1751417h} is chosen small enough. By Cauchy-Schwartz, it proves that $|T_3|\le CN^{-1}$.

Let us now treat $T_1$ and $T_2$. We have \beq
T_2&=& \E_{\ba_k}[\E_k(\eps_k)\E_k(J_{k, \diag}')]
\\
 &=& \E_{\ba_k}[\E_k(\eps_k )\E_k(J_{k, \Tr}' + EJ_{k, \diag}'J_{k, \Tr}')]
 \\
 &=& (\E_kJ_{k, \Tr}')\E_{\ba_k}\E_k(\eps_k ) +  \E_{\ba_k}[\E_k(\eps_k )\E_k (EJ_{k, \diag}'J_{k, \Tr}')]
 \\
& \leq & C\E_{\ba_k}[|\E_k(\eps_k )\E_k E|]
\eeq
where in the last two lines we used that $J_{k, \Tr}'$ is independent of $\ba_k$ and that $\E_{\ba_k}\eps_k = 0$.
By Cauchy-Schwartz and our previous $\epsilon$, we get that
\[
|T_2| \le  C\sqrt{\E_{\ba_k}|\E_k(\eps_k )|^2 \E_{\ba_k}|\E_k E|^2} \le  C N^{-\frac{1}{2}} N^{\frac{1}{2}\left(\frac{4}{\alpha} - \frac{\alpha}{4} - 1+\epsilon(4-\al)\right)}
 \ll N^{1 - \frac{\alpha}{2} }.
\]
The same bound holds for $T_1$. This concludes the proof of  Proposition \re{85141}.\epr

\subsection{Computation of the limit}\label{calculation}
By what precedes, to prove \eqre{DefC(u)513}, it suffices  to prove that the random variables  \be\la{2071215}f_k(z):=\f{1+\sum_j|\ba_k(j)|^2(G^{(k)}(z))^2_{jj}}{z-\sum_j|\ba_k(j)|^2G^{(k)}(z)_{jj} }\ee
satisfy the convergence in \pro \be\la{2071215CDD} N^{-1+\al/2} \E_{k-1}[(\E_{k}-\E_{k-1})(f_k(z))(\E_{k}-\E_{k-1})(f_k(z'))] \lto 2uC(z,z')\ee as $N,k\to\infty$ with $k/N\to u$.

 Let   $\wt{f}_k(z)$ be defined as $f_k(z)$, but with the matrix $A$ replaced by a matrix $\wt{A}$ whose entries $\wt{a}_{ij}$ are  the ones of $A$   if $i\le k$ and $j\le k$ and  independent random variables with the same distribution as the entries of $A$ if $i> k$ or $j> k$. Let $\wt{G}^{(k)}(z)$ denote the resolvent of $\wt{A}^{(k)}$. This notation is convenient as we can then express a product of integrals as integrals over different sets of variables. Let also  $\E_{\ba_k,\wt{\ba}_k}$ denote  the expectation with respect to the randomness of the $k$-th columns of $A$ and $\wt{A}$. Furthermore, $\E_k$ still denotes the conditional expectation with respect to the $\si$-algebra generated by the $k\ti k$ upper left corner of $A$ (or of $\wt{A}$, as they share the same $k\ti k$ upper left corner).

\beg{lem}\la{1951413}We have\beq && \E_{k-1}[(\E_{k}-\E_{k-1})(f_k(z))(\E_{k}-\E_{k-1})(f_k(z'))]=\\ &&\qquad\qquad\qquad\qquad\E_{k}[\E_{\ba_k,\wt{\ba}_k}(f_k(z)\wt{f}_k(z'))]-\E_{k}\E_{\ba_k}f_k(z)\E_{k}\E_{\ba_k}f_k(z').\eeq\en{lem}
\bpr
On the $\si$-algebra generated by $A$ and $\wt{A}$, we have   $\E_{k-1}=\E_{\ba_k,\wt{\ba}_k}\E_k=\E_{k}\E_{\ba_k,\wt{\ba}_k}$,  hence by Lemma \re{1912131},
 \beq && \E_{k-1}[(\E_{k}-\E_{k-1})(f_k(z))(\E_{k}-\E_{k-1})(f_k(z'))]=\\ &&\qquad\qquad\qquad\qquad\E_{\ba_k,\wt{\ba}_k}[\E_kf_k(z)\E_kf_k(z')]-\E_{k}\E_{\ba_k}f_k(z) \E_{k}\E_{\ba_k}f_k(z')\eeq
  Now, note that $\E_kf_k(z')=\E_k\wt{f}_k(z')$, hence
  \beq && \E_{k-1}[(\E_{k}-\E_{k-1})(f_k(z))(\E_{k}-\E_{k-1})(f_k(z'))]=\\ &&\qquad\qquad\qquad\qquad\E_{\ba_k,\wt{\ba}_k}[\E_kf_k(z)\E_k\wt{f}_k(z')]-\E_{k}\E_{\ba_k}f_k(z) \E_{k}\E_{\ba_k}f_k(z')\eeq
 Lastly, $$\E_kf_k(z)\E_k\wt{f}_k(z')=\E_k[f_k(z) \wt{f}_k(z')]$$ thus
 \beq && \E_{k-1}[(\E_{k}-\E_{k-1})(f_k(z))(\E_{k}-\E_{k-1})(f_k(z'))]\\ &&\qquad\qquad\qquad\qquad=\E_{\ba_k,\wt{\ba}_k}[\E_k[f_k(z) \wt{f}_k(z')]]-\E_{k}\E_{\ba_k}f_k(z) \E_{k}\E_{\ba_k}f_k(z')\\ &&\qquad\qquad\qquad\qquad=\E_k[\E_{\ba_k,\wt{\ba}_k}[f_k(z) \wt{f}_k(z')]]-\E_{k}\E_{\ba_k}f_k(z) \E_{k}\E_{\ba_k}f_k(z').\eeq
 \epr

    For $w \in \C \bck\R$,
	%\be\la{rep}
	$\ds 	\ff{w}=-\ii\op{sgn}_w   \int_{0}^{+\infty}e^{\op{sgn}_w \ii tw}\ud t,$
	 (recall that  $\mathrm{sgn}_w=\op{sgn}(\Im w)$).
	 Hence, letting $w = z-\sum_j|\ba_k(j)|^2G^{(k)}(z)_{jj}$ and using \eqre{41220151}, we get that\beq f_k(z)&=&\f{1+\sum_j|\ba_k(j)|^2(G^{(k)}(z))^2_{jj}}{z-\sum_j|\ba_k(j)|^2G^{(k)}(z)_{jj} }\\
	 &=&-\ii \, \mathrm{sgn}_z \int_{0}^{+\infty}(1+\sum_j|\ba_k(j)|^2(G^{(k)}(z)^2)_{jj})e^{\mathrm{sgn}_z\ii t(z-\sum_j|\ba_k(j)|^2G^{(k)}(z)_{jj})}\ud t\eeq
	 By \eqre{41220151} the above integral in $t$ and in the randomness of $\ba_k$ is absolutely convergent and thus we can interchange $\E_{\ba_k}$ and $\int_{t=0}^{+\infty}$. Thus
	 \beq \E_{\ba_k} f_k(z)
	 &=&-\ii \, \mathrm{sgn}_z \int_{0}^{+\infty}\E_{\ba_k}(1+\sum_j|\ba_k(j)|^2(G^{(k)}(z)^2)_{jj})e^{\mathrm{sgn}_z\ii t(z-\sum_j|\ba_k(j)|^2G^{(k)}(z)_{jj})}\ud t\eeq
Then, for any $t>0$, we have  that
	 $$\E_{\ba_k}(1+\sum_j|\ba_k(j)|^2(G^{(k)}(z)^2)_{jj})e^{\mathrm{sgn}_z\ii t(z-\sum_j|\ba_k(j)|^2G^{(k)}(z)_{jj})}=\qquad\qquad $$ $$ \qquad\qquad \E_{\ba_k}\ff{\mathrm{sgn}_z\ii t}\pa_z\{e^{\mathrm{sgn}_z\ii t(z-\sum_j|\ba_k(j)|^2G^{(k)}(z)_{jj})}\}$$ and that $\E_{\ba_k}$ and $\pa_z$ can be permuted by  \eqre{41220151} again. Hence, we have
	 \beq
\E_{\ba_k} f_k(z)  &=&-\ii \, \mathrm{sgn}_z \int_{0}^{+\infty}\ff{\mathrm{sgn}_z\ii t}\pa_z\E_{\ba_k} \{e^{\mathrm{sgn}_z\ii t(z-\sum_j|\ba_k(j)|^2G^{(k)}(z)_{jj})}\}\ud t	
	 \eeq
	 and for
$\phi_N(\lam)= \E \big[ \exp( -\ii \lambda |a_{11}|^2) \big]$
 as defined at \eqre{2071216h23},
 \be\la{50320161}
  \E_{\ba_k}f_{k} \;=\;-  \int_{0}^{+\infty} \pa_z\ff{t}e^{\mathrm{sgn}_z\ii tz}\prod_j \phi_N(\mathrm{sgn}_ztG^{(k)}(z)_{jj})\ud t.
  \ee
	 Let us now use \eqre{2071216h23} : $$  \phi_N(\lam)= 1-\f{\ii\lam\sigma_N^2}{N} + c\f{(\ii\lam \sigma_N^2)^\f{\al}{2}}{N^{\f{\al}{2}}}+\f{|\lam\sigma_N|^{\f{\al}{2}}}{ N^{\f{\al}{2}}}\eps_N(\ii\lam \sigma_N^2/N)$$ Hence  \begin{align*}\E_{\ba_k}f_{k}(z) %&=- \int_{0}^{+\infty}\pa_z \ff{t}e^{\mathrm{sgn}_z\ii tz}\prod_j \lf(1-\f{\ii \, \mathrm{sgn}_ztG^{(k)}(z)_{jj}}{N}-c\f{(\ii \, \mathrm{sgn}_ztG^{(k)}(z)_{jj})^{\f{\al}{2}}}{N^{\f{\al}{2}}}+o(N^{-\f{\al}{2}})\ri)\ud t\\
&=- \int_{0}^{+\infty}\pa_z \ff{t} e^{\mathrm{sgn}_z\ii tz}\prod_j \lf(1+\f{u_j(z)}{N/\sigma_N^2}\ri)\ud t
\end{align*}
for \beqy\la{43161} u_j(z)&:=&  (N/\sigma_N^2)\lf(\phi_N(\mathrm{sgn}_ztG^{(k)}(z)_{jj})-1\ri)
\\ \nonumber
&=&-\ii \, \mathrm{sgn}_ztG^{(k)}(z)_{jj} +c\f{(\ii\, \mathrm{sgn}_ztG^{(k)}(z)_{jj})^{\f{\al}{2}}} {(N/\sigma_N^2)^{\f{\al-2}{2}}} +
\f{|tG^{(k)}(z)_{jj}|^{\f{\al}{2}}}{(N/\sigma_N^2)^{\f{\al-2}{2}}}\eps_N(\ii \, \mathrm{sgn}_ztG^{(k)}(z)_{jj} /(N/\sigma_N^2)).\eeqy
Let \beq
\delta(z, t) &:=& \ff{t} e^{\mathrm{sgn}_z\ii tz}\prod_j \phi_N(\mathrm{sgn}_ztG^{(k)}(z)_{jj})
\\ &&\qquad -
 \ff{t}e^{\mathrm{sgn}_z\ii tz-\f{\mathrm{sgn}_z\ii t}{N/\sigma_N^2}\Tr G^{(k)}(z)}\bigg(1 + \f{c}{(N/\sigma_N^2)^{\al/2}}\sum_j (\ii\,\mathrm{sgn}_ztG^{(k)}(z)_{jj})^{\al/2}  \bigg).
\eeq
We want to show that
\[
\left|\int_0^\infty \pa_z \delta(z, t) \ud t\right| = o(N^{-\alpha/2+1}) .
\]
Since $\delta(z, t)$ is analytic in $z\in \C\bck\R$, by the Cauchy inequality
\[
|\pa_z \delta(z, t)| < \frac{2M_t}{\Im z}
\]
where $M_t = \max_{B(z, \Im z/2)} |\delta(z, t)|$. Let $z_t$ be the maximizer of $\delta(w, t)$ on $B(z, \Im z/2)$. Then
\[
\left|\int_0^\infty \pa_z \delta(z, t) \ud t\right| \leq \frac{2}{\Im z} \int_0^\infty |\delta(z_t, t)| \ud t.
\]
Let $0< \gamma <(\al-2)/(2\al)<1/2$. Then we split the above integral into two parts:
\[
\int_0^\infty |\delta(z_t, t)| \ud t = \left(\int_0^{N^\gamma} + \int_{N^\gamma}^\infty\right) |\delta(z_t, t)| \ud t.
\]
It is easy to see that
\[
\int_{N^{\gamma}}^\infty |\delta(z_t, t)| \leq Ce^{-\frac{N^{\gamma}|\Im z|}{2}}
\]
so we focus on the first integral.

Recalling that $\phi_N(\mathrm{sgn}_ztG^{(k)}(z)_{jj}) =   1+\f{u_j(z)}{N/\sigma^2_N}$, we write $\ds\int_0^{N^{\gamma}}|\delta(z_t, t)| \ud t$ as the integral of a sum of three errors $\delta_1, \delta_2, \del_3$:
\begin{equation*}\begin{split}
&\int_0^{N^{\gamma}} |\delta(z_t, t) |\ud t
\leq
 \int_0^{N^{\gamma}} \left| \ff{t} e^{\mathrm{sgn}_z\ii tz_t}\left(\prod_j \lf(1+\f{u_j(z_t)}{N/\sigma^2_N}\ri)- e^{\ff{N/\sigma_N^2}\sum_j u_j(z_t)}\right)\right|
 \\
 &
+\int_0^{N^{\gamma}} \bigg| \ff{t} e^{\mathrm{sgn}_z\ii tz_t}\bigg(e^{\ff{N/\sigma_N^2}\sum_j u_j(z_t)} \\
&\qquad\qquad\qquad\qquad - e^{-\f{\mathrm{sgn}_z\ii t}{N/\sigma_N^2}\Tr G^{(k)}(z_t)}\times\\
&\qquad\qquad\qquad\qquad\bigg(1 +\f{c}{(N/\sigma_N^2)^{\al/2}}\sum_j (\ii\,\mathrm{sgn}_ztG^{(k)}(z_t)_{jj})^{\al/2}  \\
&\qquad\qquad\qquad\qquad+ \sum_j\f{|tG^{(k)}(z_t)_{jj}|^{\f{\al}{2}}}{(N/\sigma_N^2)^{\f{\al}{2}}}\eps_N(\ii \, \mathrm{sgn}_ztG^{(k)}(z_t)_{jj} /(N/\sigma_N^2)) \bigg)
\bigg)\bigg|
\\
 &
+\int_0^{N^{\gamma}} \left| \ff{t} e^{\mathrm{sgn}_z\ii tz}e^{-\f{\mathrm{sgn}_z\ii t}{N/\sigma_N^2}\Tr G^{(k)}(z_t)}
\left(
\sum_j\f{|tG^{(k)}(z_t)_{jj}|^{\f{\al}{2}}}{(N/\sigma_N^2) ^{\f{\al}{2}}}\eps_N(\ii \, \mathrm{sgn}_ztG^{(k)}(z_t)_{jj} /(N/\sigma_N^2)) \right)\right|
\\
& =: \delta_1 + \delta_2 + \delta_3
\end{split}\end{equation*}
To get a bound on $\delta_1$ we use Lemma \ref{2541400h15} for each $t$ with
\begin{equation*}\begin{split}
M_t &:= \max_i |u_i(z_t)|\sigma_N^2
\\
&
 = \max_i |-\ii \, \mathrm{sgn}_ztG^{(k)}(z_t)_{jj} +c\f{(\ii\, \mathrm{sgn}_ztG^{(k)}(z_t)_{jj})^{ \f{\al}{2}}}{(N/\sigma_N^2)^{\f{\al-2}{2}}}+\f{|tG^{(k)}(z_t )_{jj} |^{\f{\al}{2}}}{(N/\sigma_N^2)^{\f{\al-2}{2}}}\eps_N(\ii \, \mathrm{sgn}_ztG^{(k)}(z_t)_{jj} /(N/\sigma_N^2))|\sigma_N^2
 \\
 &
 \leq Ct
\end{split}\end{equation*}
for a constant $C$.  Besides, $|\phi_N(\mathrm{sgn}_ztG^{(k)}(z)_{jj})|\le 1$, hence  by \eqre{43161}, $\Re(u_j(z_t)\sigma_N^2)\le 0$.  Hence
\[\delta_1 \leq \int_0^{N^\gamma}\left| \frac{1}{t}e^{\mathrm{sgn}_z\ii tz_t}\f{C^2t^2}{N}e^{\ff{N}\sum_j \Re(u_j(z_t)\sigma_N^2)+ \f{C^2t^2}{N}} \right|\ud t \leq \f{C^2}{N}\int_0^{N^\gamma}\left|te^{\mathrm{sgn}_z\ii tz_t} \ri|e^{\f{C^2t^2}{N}} \ud t \leq \frac{C}{N^{1-2\gamma}},
\] where we used the fact that $\ga<1/2$.

To get a bound on $\delta_2$, we can use the Taylor series expansion to check that for $ |x| \le  1/2$
\[
|e^{x} - (1+x)| = \left|x^2\sum_{k = 0}^\infty \frac{x^k}{(k+2)!}\right| \le  |x|^2\sum_{k = 0}^\infty \frac{|x|^k}{(k+2)!}\le \frac{|x|^2}{1-|x|}\le  2|x|^2
\]
This yields that
\begin{equation*}\begin{split}
\exp&\left(c\sum_j\f{(\ii\, \mathrm{sgn}_ztG^{(k)}(z_t)_{jj})^{ \f{\al}{2}}}{(N/\sigma_N^2)^{\alpha/2}}+ \sum_j\f{|tG^{(k)}(z_t )_{jj} |^{\f{\al}{2}}}{(N/\sigma_N^2)^{\alpha/2}}\eps_N(\ii \, \mathrm{sgn}_ztG^{(k)}(z_t)_{jj} /(N/\sigma_N^2))\right)\\
&-
 \left(1+c\sum_j \f{(\ii\, \mathrm{sgn}_ztG^{(k)}(z_t)_{jj})^{ \f{\al}{2}}}{(N/\sigma_N^2)^{\alpha/2}}+\sum_j\f{|tG^{(k)}(z_t )_{jj} |^{\f{\al}{2}}}{(N/\sigma_N^2)^{\alpha/2}}\eps_N(\ii \, \mathrm{sgn}_ztG^{(k)}(z_t)_{jj} /(N/\sigma_N^2))\right)
 \\
 &
\le
 2 \left| c\sum_j\f{(\ii\, \mathrm{sgn}_ztG^{(k)}(z_t)_{jj})^{ \f{\al}{2}}}{(N/\sigma_N^2)^{\alpha/2}}+ \sum_j\f{|tG^{(k)}(z_t )_{jj} |^{\f{\al}{2}}}{(N/\sigma_N^2)^{\alpha/2}}\eps_N(\ii \, \mathrm{sgn}_ztG^{(k)}(z_t)_{jj} /(N/\sigma_N^2)) \right|^2
\end{split}\end{equation*}
so that, as $\ga<(\al-2)/(2\al)$ implies that $\gamma\alpha - \alpha+2< -\alpha/2+1$ and for some constants $C_1, C_2$,
\begin{equation*}\begin{split}
\delta_2 &\leq C_1\int_0^{N^{\gamma}}  \ff{t} e^{-\f{\mathrm{sgn}_z\ii t}{N}\Tr G^{(k)}(z_t)}\ti \\ & \qquad\qquad
\left| c\sum_j \f{(\ii\, \mathrm{sgn}_ztG^{(k)}(z_t)_{jj})^{ \f{\al}{2}}}{N^{\alpha/2}}
+
\sum_j\f{|tG^{(k)}(z_t )_{jj} |^{\f{\al}{2}}}{N^{\alpha/2}}\eps_N(\ii \, \mathrm{sgn}_ztG^{(k)}(z_t)_{jj} /N) \right|^2\ud t\\
\\
&\leq C_2 N^{\gamma\alpha - \alpha+2} \ll N^{-\alpha/2+1}
\end{split}\end{equation*}

To get a bound on $\delta_3$ we do a dyadic decomposition of the integral. We integrate on $[2^k, 2^{k+1}]$ with $k$ such that $2^{k}< N^{\gamma}$.
This yelds
\begin{equation*}\begin{split}
\int_{2^k}^{2^{k+1}} &\left| \ff{t} e^{\mathrm{sgn}_z\ii tz}e^{-\f{\mathrm{sgn}_w\ii t}{(N/\sigma_N^2)}\Tr G^{(k)}(z_t)}\left(\sum_j\f{|tG^{(k)}(z_t)_{jj}|^{\f{\al}{2}}} {(N/\sigma_N^2)^{\f{\al}{2}}}\eps_N(\ii \, \mathrm{sgn}_ztG^{(k)}(z_t)_{jj} /(N/\sigma_N^2)) \right)\right|dt \\
&
\leq C 2^ke^{-|C\Im z|2^k}\frac{(2^k)^{\alpha/2-1}}{N^{\alpha/2-1}} \max_{t \in [0, N^\gamma]}\eps_N(\ii \, \mathrm{sgn}_ztG^{(k)}(z_t)_{jj} /(N/\sigma_N^2))
\end{split}\end{equation*}
Noting that $\sum_k 2^{k\alpha/2}e^{-|C\Im z|2^k}$ is convergent we get that
\[
\delta_3\leq C \frac{\max_{t\in [0, N^\gamma]}\eps_N(\ii \, \mathrm{sgn}_ztG^{(k)}(z_t)_{jj} /(N/\sigma_N^2))}{N^{\alpha/2 - 1}} \ll N^{-\alpha/2 + 1}
\]

As a consequence,
$$\E_{\ba_k}f_{k}(z) =- \int_{0}^{+\infty} \pa_z\ff{t}e^{\mathrm{sgn}_z\ii tz-\f{\mathrm{sgn}_z\ii t}{N/\sigma_N^2}\Tr G^{(k)}(z)}\bigg(1+\f{c}{(N/\sigma_N^2)^{\al/2}}\sum_j (\ii\,\mathrm{sgn}_ztG^{(k)}(z)_{jj})^{\al/2} \bigg) \ud t +o(N^{-\al/2+1})
$$

Then Lemma \re{lemSCdiagGreen}, whose statement and proof are postponed until the next section,  implies that the diagonal terms $\sigma_N^2G^{(k)}(z)_{jj}$ in the previous expression are close to the Stieltjes transform $\Gsc(z)$ of the semicircle law with support $[-2,2]$. It also implies that $N^{-1}\sigma_N^2\Tr G^{(k)}(z)$ is close to  $\Gsc(z)$. As before, the factor $e^{\mathrm{sgn}_z\ii tz}$ allows us to take the integral $\ud t$ to infinity.
It follows   that \be\la{195141}\E_{\ba_k}f_k(z)=- \int_{0}^{+\infty} \pa_z\ff{t}e^{\mathrm{sgn}_z\ii tz-\mathrm{sgn}_z\ii t\Gsc(z)}\bigg(1+ N^{-\al/2+1}cK(z,t)^{\al/2}\bigg) \ud t+o(N^{-\al/2+1}) ,\ee
where $o(1)$ is for the convergence in probability.

Of course, the same holds for $\E_{\ba_k}f_k(z')$.

Let us now compute $\E_{\ba_k,\wt{\ba}_k}(f_k(z)\wt{f}_k(z'))$.
First, as above,  $$f_k(z) \wt{f}_k(z')= \int_{t,t'>0}\pa_z\pa_{z'}\ff{tt'} \{e^{\mathrm{sgn}_z\ii t(z-\sum_j|\ba_k(j)|^2G^{(k)}(z)_{jj})+\mathrm{sgn}_{z'}\ii t'(z'-\sum_j|\wt{\ba}_k(j)|^2\wt{G}(z')^{(k)}_{jj})}\}\ud t	 \ud t'
$$
	 Hence as $\E_{\ba_k}$ denotes the integration with respect to the $k$-th columns of $A$ and $A'$ and as these columns are identical up to the $k$-th entry and independent from the $k+1$-th entry on, we have
\begin{align*}\E_{\ba_k,\wt{\ba}_k}(f_k(z)\wt{f}_k(z'))&=\int_{t,t'>0}\pa_z\pa_{z'}\ff{tt'}\{e^{\mathrm{sgn}_z\ii tz+\mathrm{sgn}_{z'}\ii t'z'}\prod_{j<k} \phi_N(\mathrm{sgn}_ztG^{(k)}(z)_{jj}+\mathrm{sgn}_{z'}t'\wt{G}(z')^{(k)}_{jj})\\ &\qquad\prod_{j>k} \phi_N(\mathrm{sgn}_ztG^{(k)}(z)_{jj})\phi_N(\mathrm{sgn}_{z'}t'\wt{G}(z')^{(k)}_{jj})\}\ud t \ud t'
%\\
%&=\int_{t,t'>0}\pa_z\pa_{z'}\ff{tt'} e^{\mathrm{sgn}_z\ii tz-\f{\mathrm{sgn}_z\ii t}{N}\Tr G^{(k)}(z)}e^{\mathrm{sgn}_{z'}\ii t'z'-\f{\mathrm{sgn}_{z'}\ii t'}{N}\Tr \wt{G}(z')^{(k)}}\\ &\qquad
%\bigg(1-\f{c}{N^{\al/2}}\sum_{j<k} (\ii\,\mathrm{sgn}_ztG^{(k)}(z)_{jj}+\ii\,\mathrm{sgn}_{z'}t'\wt{G}(z')^{(k)}_{jj})^{\al/2}  \\ &-\f{c}{N^{\al/2}}\sum_{j>k} (\ii\,\mathrm{sgn}_ztG^{(k)}(z)_{jj})^{\al/2}+(\ii\,\mathrm{sgn}_{z'}t'\wt{G}(z')^{(k)}_{jj})^{\al/2} +o(N^{-\al/2+1}) \bigg) \ud t \ud t'
\end{align*}
Then, proceeding as above when we computed $\E_{\ba_k} f_k(z)$ (from   \eqre{50320161} to   \eqre{195141}), we get that
\beq\la{1951412}\E_{\ba_k,\wt{\ba}_k}(f_k(z)\wt{f}_k(z')) &=&\int_{t,t'>0}\pa_z\pa_{z'}\ff{tt'} e^{\mathrm{sgn}_z\ii tz-K(z,t)}e^{\mathrm{sgn}_{z'}\ii t'z'-K(z',t')}
\\&& \nonumber
\bigg(1+cuN^{-\al/2+1} (K(z,t)+K(z',t'))^{\al/2}  \\ \nonumber
&&+c(1-u)N^{-\al/2+1} (K(z,t)^{\al/2}+K(z',t')^{\al/2})  +o(N^{-\al/2+1}) \bigg) \ud t \ud t'
\eeq

 This equation, together with \eqre{195141} and Lemma \re{1951413}, imply \eqre{2071215CDD}. This concludes the proof.

% It implies that \beq &&\E_{k}[\E_{\ba_k}(f_k\ti f_k'')]-\E_{k}\E_{\ba_k}f_k\ti \E_{k}\E_{\ba_k}f_k=o(N^{-\al/2+1}) -\\ \\ &&\int_{t,t'>0}\ud t\ud t'\pa_z\pa_{z'}\ff{tt'} e^{\mathrm{sgn}_z\ii tz-\f{\mathrm{sgn}_z\ii t}{N}\Tr G^{(k)}(z)} e^{\mathrm{sgn}_{z'}\ii t'z'-\f{\mathrm{sgn}_{z'}\ii t'}{N}\Tr G'_k(z')}  \f{c}{N^{\al/2}}S_N(k,z,t,t')
%\eeq
%where $$S_N(k,z,t,t'):=$$ $$\sum_{j<k} \{(\ii\,\mathrm{sgn}_ztG^{(k)}(z)_{jj}+i\mathrm{sgn}_{z'}t'\wt{G}(z')^{(k)}_{jj})^{\al/2}
%- [(\ii\,\mathrm{sgn}_ztG^{(k)}(z)_{jj})^{\al/2}+(i\mathrm{sgn}_{z'}t'\wt{G}(z')^{(k)}_{jj})^{\al/2}]$$

\subsection{Concentration of the diagonal terms of the resolvent}\label{concentration}
 %We shall here use the following lemma, which states, among other things, the convergence of the empirical spectral law of $A$ towards the semicircle law (that was already stated, see e.g. \cite{bai-silver-book}).

\beg{lem}\la{lemSCdiagGreen}
For any fixed $z\in \C\bck\R$, for any $p\ge 1$, and any $j\in \{1, \ld, N\}$,  the sequence $$ \|G(z)_{jj}-\Gsc(z)  \|_{L^p},$$ tends to zero.
\en{lem}

\begin{proof}  As $|G(z)_{jj}-\Gsc(z)|\le 2(\Im z)^{-1}$, it suffices to prove the result for any value of $p$. By the Schur complement formula (see \cite[Th. 11.4]{bai-silver-book}), we know that \be\la{MumfordSons} G(z)_{jj}=\ff{z-\f{x_{jj}}{\sqrt{N}}-\ff{N}\sum_{i,k : i\ne j, k\ne j}G^{(j)}(z)_{ik}x_{ji}x_{kj}},\ee where $G^{(j)}(z)=(z-A^{(j)})^{-1}$ and $A^{(j)}$ is the matrix obtained after removing the $j$th row and the $j$th column of $A$.

We recall that the denominator is equal to
\[ z - \ff{N}\Tr G(z) + \ff{N}(\Tr G^{(j)}(z)- \Tr G(z)) + \eta_j + E ,  \]
where $\eta_j$ is as in \eqref{d:etak} and $E = \frac{1}{N}\sum G_{jj}^{(j)} \left(|x_{jj}|^2 - 1\right)$.

We will show that for a certain choice of $0<p<1$ and a certain choice of $\epsilon_0 >0$,
\begin{equation}\label{e:EEnu}
\E |E|^p = o(N^{-\epsilon_0}).
\end{equation}

Recall that here $x_{jj}$ has been truncated at $N^{\frac{1}{\alpha} + \frac{1}{4} + \epsilon}$. In the proof of this Lemma we will truncate the variables further at $N^{\beta_1}.$ Since each variable is truncated, independence of variables is retained. We let $\Omega_k$ be event that for $k$ of $j$'s, we have that $|x_{jj}| \geq N^{\frac{1}{4}-\delta}$. Then, by Jensen's inequality and Lemma \ref{2151411h50},
\begin{equation}
\E_{\Omega_0}\left|\frac{1}{N}\sum G_{jj}(|x_{jj}|^2-1)\right|^{p}
\leq \left(\E_{\Omega_0}|E|^2\right)^{p/2} \leq C(N^{\beta_1(4-\al)-1})^{p/2}.
\end{equation}
Now,
\begin{multline}
E_{\Omega_0^C}\left|\frac{1}{N}\sum G_{jj}(|x_{jj}|^2-1)\right|^{p} \leq CN \mathbb{P}(|x_{11}| >N^{\beta_1})N^{\left(\frac{1}{4}+ \frac{1}{\alpha} + \epsilon\right)p}\leq CN^{-\alpha \beta_1 + \left(\frac{1}{4}+ \frac{1}{\alpha} + \epsilon\right)p}\\
\end{multline}
Choosing $p$ and $\beta_1$ in such a way that $\beta_1(4-\al)-1<0$ and $-\alpha \beta_1 + \left(\frac{1}{4}+ \frac{1}{\alpha} + \epsilon\right)p< 0$ finishes the proof of \eqref{e:EEnu}.

Again, using Jensen's Inequality we obtain that
\[\E|\eta_j|^p \leq (\E|\eta_j|^2)^{p/2} \leq (CN)^{-p/2}.
\]
Lastly, using that $| \Tr G^{(j)}(z)- \Tr G(z)|\le C$ (see Lemma \re{lat}),  the denominator of the RHS of \eqre{MumfordSons} can be written $$z-\ff{N}\Tr G(z)+ O_{L^p}(N^{-\epsilon_0})$$.

As the function $f_z:x\mapsto\ff{z-x}$ has uniformly bounded gradient on the half plane $\{x\in \C\ste \Im z \Im x >0\}$, we get that $$G(z)_{jj}=\ff{z-\ff{N}\Tr G(z)}+O_{L^p}(N^{-\epsilon_0}).$$
    Besides, by  \cite[Lem. C.1]{alice-charles-HT}, we know that $$\ff{N}\Tr G(z)-\E[\ff{N}\Tr G(z)]= O_{L^p}(N^{-\epsilon_0}),$$ hence
$$G(z)_{jj}=\ff{z-\E[\ff{N}\Tr G(z)]}+ O_{L^p}(N^{-\epsilon_0}).$$
At last, by \cite[Th. 2.5]{bai-silver-book}, we know that for any fixed $z\in \C\bck \R$, $$\E[\ff{N}\Tr G(z)]-\Gsc(z)\lto 0.$$
We conclude by using the fact that $\ds \ff{z-\Gsc(z)}=\Gsc(z)$.\epr

 \section{Appendix}\subsection{Quadratic forms in heavy-tailed variables}
 \beg{lem}\la{2151411h50}Let $\ba=(a_1, \ld, a_N)^T$ be a column vector whose entries are i.i.d., centered and satisfy (ii) and (iii) of Lemma \re{215141}. Then for any deterministic matrix $G$, the random variables $$X:=\sum_{i\ne j} G_{ij} \ovl{a_i} a_j \qquad\qquad E :=\sum_{i } G_{ii} |a_{i}|^2-\ff{N}\Tr G$$ satisfy $$\E[|X|^2] \le  2N^{-2}\Tr(GG^*)\le 2N^{-1}\|G\|^2\qquad\qquad \E[|E|^2]\le 10C(\|G\|^2+1)N^{\beta(4-\al)-1}.$$
 \en{lem}

 \bpr Direct computations, the second one using that $a_j$ is centred for the first one and that $|a_j|^2 - \frac{1}{N}$ is centred for the second one.
 \epr

 \beg{rmk}We shall sometimes use this lemma after removal of the $k$th row and column of $B$  and of the $k$th entry of $\ba$, but it suffices to apply the lemma with the matrix deduced from $B$ by setting its $k$th row and column to zero.
 \en{rmk}

 \subsection{CLT for martingales}

Let $(\mc{F}_k )_{k\ge 0}$ be a filtration \st $\mc{F}_0 =\{\emptyset,\Omega\}$ and let $(M_k )_{k\ge 0}$ be a square-integrable complex-valued martingale starting at zero with respect to this filtration. For $k\ge 1$, we define the random variables $$Y_k:=M_k-M_{k-1}\qquad v_k:=\E[|Y_k|^2\,|\,\mc{F}_{k-1} ]\qquad  \tau_k:=\E[Y_k^2\,|\,\mc{F}_{k-1} ]$$ and we also define $$v:=\sum_{k\ge 1}v_k\qquad  \tau:=\sum_{k\ge 1} \tau_k\qquad L(\eps):=\sum_{k\ge 1} \E[|Y_k|^2\one_{|Y_k|\ge \eps}].$$

Let now everything depend on a parameter $N$, so that $\mc{F}_k=\mc{F}_k(N), M_k=M_k(N), Y_k=Y_k(N),v=v(N),\tau=\tau(N),  L(\eps)=L(\eps, N), \ldots$

Then we have the following theorem. It is proved in the real case at \cite[Th. 35.12]{Billingsley}.  The complex  case can be deduced noticing that for $z\in \C$, $\Re(z)^2, \Im(z)^2$ and $\Re(z)\Im(z)$ are linear combinations of $z^2$, $\ovl{z}^2$, $|z|^2$.
\beg{Th}\la{thconvmart}Suppose that for  some  constants $v\ge  0, \tau\in \C$, we have the convergence in probability $$v(N)\Ninf v\qquad \tau(N)\Ninf \tau$$ and that for each $\eps>0$, $$L(\eps, N)\Ninf 0.$$ Then we have the convergence in distribution $$M_N(N)\Ninf Z,$$ where $Z$ is a centered complex Gaussian variable \st $\E(|Z|^2)=v$ and $ \E(Z^2)=\tau$.
\en{Th}

To apply this theorem, we shall use the following lemma.
\beg{lem}\la{1912131}Let $\E_{k-1}$ and $\E_k$ denote conditional expectations given some $\si$-algebras $\F_{k-1}\subset\F_k$. then for any pair $A,B$ of $L^2$ r.v., we have $$\E_{k-1}[(\E_k-\E_{k-1})(A)(\E_k-\E_{k-1})(B)]=\E_{k-1}[\E_k(A)\E_k(B)]-\E_{k-1}(A)\E_{k-1}(B).$$
\en{lem}

\bpr  We have %using the fact that $\E_{k-1}\E_k=\E_k\E_{k-1}=\E_{k-1}\E_{k-1}=\E_{k-1}$,
\beq\E_{k-1}[(\E_k-\E_{k-1})(A)(\E_k-\E_{k-1})(B)]&=&\E_{k-1}[\E_kA (\E_k-\E_{k-1})(B)]\\ &&
-\E_{k-1}[\E_{k-1}(A) (\E_k-\E_{k-1})(B)]\\
&=&\E_{k-1}[\E_kA \E_k B]\\ &&-\underbrace{\E_{k-1}[\E_kA \E_{k-1}(B)]}_{=\E_{k-1}(B)\E_{k-1}[\E_kA]=\E_{k-1}(A)\E_{k-1}(B)}\\ &&
-\E_{k-1}(A) \underbrace{\E_{k-1}[(\E_k-\E_{k-1})B]}_{=0}
\eeq
which concludes the proof.
\epr

\subsection{A lemma about large products and the exponential function}
\beg{lem}\la{2541400h15} Let $u_i$, $i=1, \ld, N$, be some complex numbers and set
 $$P:=\prod_{i=1}^N (1+\f{u_i}{N})\qquad\qquad S:=\ff{N}\sum_{i=1}^Nu_i\qquad\qquad M:=\max_i|u_i|.$$ There is a universal constant $R>0$ (independent of $N$ and of $M$) \st   $$  \f{M}{N}\le R\;\implies\; |P-e^S|\le \f{M^2}{N}e^{\Re (S)+\f{M^2}{N}} .$$
\en{lem}

\bpr Let $L(z)$ be defined on $B(0,1)$ by $\log(1+z)=z+z^2L(z)$ and $R>0$ be \st on $B(0,R)$, $|L(z)|\le 1.$
If $\f{M}{N}\le R$, we have  $$P\ =\ \prod_i\exp\lf\{\f{u_i}{N}+\f{u_i^2}{N^2} L\left(\f{u_i}{N}\right)\ri\}\
=\ e^S\exp\lf\{\sum_i \f{u_i^2}{N^2} L\left(\f{u_i}{N}\right) \ri\},$$so that $$P-e^S\ =\ e^S\lf(\exp\lf\{\sum_i \f{u_i^2}{N^2} L\left(\f{u_i}{N}\right) \ri\}-1\ri)$$
Since for any $z$, $|e^{z}-1|\le |z|e^{|z|},$ the conclusion follows. \epr

\subsection{Linear algebra}
Let $H=[h_{ij}]$ be an $N\ti N$ Hermitian matrix and  $z\in \C\bck\R$. Define $G:=(z-H)^{-1}$.

\beg{lem}[Difference of traces of a matrix and its major submatrices]\la{lat} Let $H_k$ be the submatrix of $H$ obtained by removing its $k$-th row and $k$-th column and set $G_k:=(z-H_k)^{-1}$. Let also $\ba_k$ be the $k$-th column of $H$ where the $k$-th entry has been removed. Then \be\la{307121}\Tr(G)-\Tr(G_k)=\f{1+\ba_k^*G_k^{2}\ba_k}{z-h_{kk}-\ba_k^*G_k\ba_k}.\ee Moreover, \be\la{307122}|\Tr(G)-\Tr(G_k)|\le \pi|\Im z|^{-1}.\ee \en{lem}

%\beg{lem}\la{lem267121} With the notation introduced above the previous lemma, for each $1\le j\le N$,  \be\la{297122}\Im z \ti \Im G_{jj}<0,\ee  \be\la{297124}  |\Im z|\ti |(G^2)_{jj}| \le |\Im G_{jj}|\le |\Im z|^{-1} \ee
%and for any $\ba=(\ba_1,\ld, \ba_N)\in \C^N$,
%\be\la{297125}
%\lf|\f{1+\sum_j|\ba_j|^2(G^2)_{jj}}{z-\sum_j|\ba_j|^2G_{jj} }\ri|\le 2|\Im z|^{-1}.\ee
%\en{lem}

\end{document}